\documentclass[12pt,a4paper]{amsart}
\usepackage{amssymb,amsmath,amsthm}
\usepackage{graphicx}
\usepackage{hyperref}
 \usepackage[color=green!40]{todonotes}

\newcommand\cA{{\mathcal A}}
\newcommand\cB{{\mathcal B}}
\newcommand\cC{{\mathcal C}}

\newcommand\cF{{\mathcal F}}

\newcommand\cG{{\mathcal G}}
\newcommand\cH{{\mathcal H}}

\makeatletter
\newtheorem*{rep@theorem}{\rep@title}
\newcommand{\newreptheorem}[2]{%
\newenvironment{rep#1}[1]{%
 \def\rep@title{#2 \ref{##1}}%
 \begin{rep@theorem}}%
 {\end{rep@theorem}}}
\makeatother

\theoremstyle{plain}
\newtheorem{theorem}{Theorem}[section]
\newreptheorem{theorem}{Theorem}
\newtheorem{lemma}[theorem]{Lemma}

\newtheorem{conjecture}[theorem]{Conjecture}

\newtheorem{observation}[theorem]{Observation}

\theoremstyle{definition}

\newtheorem{defn}[theorem]{Definition}

\newtheorem{fact}[theorem]{Fact}

\newcommand\cref[1]{Corollary~\ref{cor:#1}}

\title{On the sizes of $t$-intersecting $k$-chain-free families}

\author{J\'ozsef Balogh}
\address{University of Illinois at Urbana-Champaign}
\email{jobal@illinois.edu}
\thanks{Balogh's research is partially
supported by NSF grants DMS-1764123 and RTG DMS-1937241, the Arnold O. Beckman Research Award (UIUC
Campus Research Board RB 18132), the Langan Scholar Fund (UIUC), and the Simons Fellowship.}

\author{William B. Linz}
\address{University of Illinois at Urbana-Champaign}
\email{wlinz2@illinois.edu}
\thanks{Linz's research is partially supported by 
 RTG DMS-1937241.}

\author{Bal\'azs Patk\'os}
\address{Alfr\'ed R\'enyi Institute of Mathematics}
\email{patkos@renyi.hu}
\thanks{Patk\'os's research is partially supported by NKFIH grants SNN 129364 and FK 132060.}

\date{}

\begin{document}

\maketitle

\begin{abstract}
    A set system $\cF$ is $t$-{\it intersecting}, if the size of the intersection of every pair of its elements has size at least $t$. A set system $\cF$ is $k$-{\it Sperner}, if it does not contain a chain of length $k+1$.
  
    Our main result is the following: Suppose that $k$ and $t$ are fixed positive integers, where  $n+t$ is even with $t\le n$ and $n$ is large enough. If $\cF\subseteq 2^{[n]}$ is a $t$-intersecting $k$-Sperner family, then $|\cF|$ has size at most the size of the sum of $k$ layers, of sizes $(n+t)/2,\ldots, (n+t)/2+k-1$. 
    This bound is best possible. The case when    $n+t$ is odd remains open.
\end{abstract}

\section{Introduction}

\subsection{Definitions and Notation}

For a positive integer $n$, we write  $[n]:= \{1, 2, \ldots, n\}$ and $2^{[n]}$ for the power set of $[n]$. For a set $S$, we denote by $\binom{S}{i}$ the family of all $i$ element subsets of $S$.

For a family of sets $\cF \subseteq 2^{[n]}$, we define $\cF_i:= \{F\in \cF: |F| = i\}$ and $f_i:= |\cF_i|$. We use $\Delta_i$ and $\nabla_i$ to denote the $i$-shadow and $i$-shade of $\cF$, respectively, so that $\Delta_i\cF := \{A: |A|=i, \hspace{2mm} A\subset F \text{ for some } F\in \cF\}$ and $\nabla_i\cF := \{A: |A|=i, \hspace{2mm} A\supset F \text{ for some } F\in \cF\}$. If the subscript $i$ is unspecified, then assuming $\cF$ is $r$-uniform, $\Delta\cF = \Delta_{r-1}\cF$ and similarly $\nabla\cF = \nabla_{r+1}\cF$.

\begin{defn}\label{defn:kspernerdefn}[$k$-Sperner family]\\
A \emph{$(k+1)$-chain} is a collection of $k+1$ sets $A_0, A_1, \ldots, A_k$ such that $A_0 \subset A_1 \subset \ldots \subset A_k$. A family of sets $\cF\subseteq 2^{[n]}$ is a \emph{$k$-Sperner family} if there is no $(k+1)$-chain  in $\cF$. If $k=1$, then $\cF$ is simply called a \emph{Sperner family} or an \emph{antichain}. 
\end{defn}

\begin{defn}\label{defn:tintdefn}[$t$-intersecting family]\\
A family of sets $\cF \subseteq 2^{[n]}$ is \emph{$t$-intersecting} if for every pair of sets $A, B \in \cF$, we have $|A\cap B| \ge t$. If $t=1$, then we write that $\cF$ is \emph{intersecting}. 
\end{defn}

\subsection{History}
The maximum size of an antichain in $2^{[n]}$ was determined by Sperner~\cite{S}. 

\begin{theorem}[Sperner]\label{spernerthm}
Let $\cF \subseteq 2^{[n]}$ be an antichain. Then, 
\[|\cF| \le \binom{n}{\lfloor{\frac{n}{2}\rfloor}}.\]
Furthermore, equality holds only if $\cF$ is one of the largest layers in the Boolean lattice $2^{[n]}$. 
\end{theorem}

Sperner's theorem was extended to $k$-Sperner families by Erd\H{o}s~\cite{E}. 

\begin{theorem}[Erd\H{o}s]\label{erdosthm}
The maximum-size $k$-Sperner family $\cF \subseteq 2^{[n]}$ is the union of the largest $k$ layers in the Boolean lattice $2^{[n]}$. 
\end{theorem}

A different extension of Sperner's theorem was given by Milner~\cite{M}. Milner additionally required the family $\cF$ to be $t$-intersecting. 

\begin{theorem}[Milner]\label{milnerthm}
If $\cF \subseteq 2^{[n]}$ is a $t$-intersecting antichain, then 
\[|\cF| \le \binom{n}{\lfloor{\frac{n+t+1}{2}\rfloor}}.\]
\end{theorem}

In a different direction, Frankl~\cite{F1} determined the maximum size of an intersecting $k$-Sperner family. Different proofs were given by Gerbner \cite{G} and by Gerbner, Methuku and Tompkins~\cite{GMT}. 

\begin{theorem}[Frankl]
Let $\cF \subseteq 2^{[n]}$ be an intersecting, $k$-Sperner family. Then, 
\[|\cF| \le \begin{cases} \sum_{i=\frac{n+1}{2}}^{\frac{n+1}{2}+k-1}\binom{n}{i}, & \text{if $n$ is odd, } \\
\binom{n-1}{\frac{n}{2}-1}+\sum_{i=\frac{n}{2}+1}^{\frac{n}{2}+k-1}\binom{n}{i}+\binom{n-1}{\frac{n}{2}+k},& \text{if $n$ is even. }\\
\end{cases}\]
Furthermore, if $n$ is odd, equality holds only if \[ \cF = \binom{[n]}{\lfloor{\frac{n}{2}\rfloor}+1} \cup \binom{[n]}{\lfloor{\frac{n}{2}\rfloor}+2} \cup \ldots \cup \binom{[n]}{\lfloor{\frac{n}{2}\rfloor}+k},\] while if $n$ is even and $k > 1$, equality holds only if for some $x \in [n]$, 
\[\cF = \left\{F \in \binom{[n]}{\frac{n}{2}}: x\in F\right\} \cup \binom{[n]}{\frac{n}{2}+1}\cup \ldots \cup \binom{[n]}{\frac{n}{2}+k-1} \cup \left\{F\in \binom{[n]}{\frac{n}{2} + k}: x\notin F\right\}.\]
\end{theorem}

A common generalization of the theorems of Milner and Frankl would be to determine the maximum size of a $t$-intersecting, $k$-Sperner family. 

Frankl~\cite{F} proposed conjectures on  the maximum size of a $t$-intersecting $k$-Sperner family $\cF\subset 2^{[n]}$ and made some progress towards proving these conjectures.  The conjectured extremal family depends on the parity of $n+t$.

In the case when $n+t$ is even, the conjectured maximum size of a $t$-intersecting, $k$-Sperner family is very easy to describe. 

\begin{conjecture}[Frankl]\label{evencase}
If $n+t$ is even, $n > t$, and $\cF \subseteq 2^{[n]}$ is a $t$-intersecting, $k$-Sperner family, then 

\[|\cF| \le \sum_{i=0}^{k-1}\binom{n}{\frac{n+t}{2}+i}.\]

\end{conjecture}
Conjecture~\ref{evencase} is clearly tight if true, as evidenced by the family $\bigcup_{i=0}^{k-1}\binom{[n]}{\frac{n+t}{2}+i}$.

The conjectured extremal families do not have such a simple structure when $n+t$ is odd. We construct two plausible candidates for the maximum size $t$-intersecting, $k$-Sperner family: 
\[\cA(t, k) = \left\{F \in \binom{[n]}{\frac{n+t-1}{2}}: n\notin F\right\} \cup \left\{A: \frac{n+t-1}{2}+1 \le |A| \le \frac{n+t-1}{2}+(k-1)\right\}.\]

\[\cB(t, k) = \left\{F \in \binom{[n]}{\frac{n+t-1}{2}}: [1, t] \in F\right\} \cup \left\{A: \frac{n+t-1}{2}+1 \le |A| \le \frac{n+t-1}{2}+(k-1)\right\} \]
\[\cup\left(\left\{B: |B| = \frac{n+t-1}{2} + k\right\} \setminus \left\{B: |B|=\frac{n+t-1}{2} + k, [1,t] \in B\right\}\right).\]

It is not hard to show that $|\cB(t,k)| \gg |\cA(t, k)|$ for $n$ sufficiently large (in terms of $k$ and $t$). However, it may be checked by computer that $\cA(t,k)$ is optimal for small values of $n$ and specific choices of $t$ and $k$, for example $t=2$ and $k=2$. We conjecture that $\cB(t, k)$ is the largest such family when $n$ is sufficiently large. 

\begin{conjecture}\label{oddcase}
There exists a positive integer $n_0=n_0(k, t)$ such that if $n+t$ is odd, $n > n_0$, and $\cF \subseteq 2^{[n]}$ is a $t$-intersecting, $k$-Sperner family, then 

\[|\cF| \le |\cB(t, k)| =  \binom{n-t}{\frac{n-t-1}{2}} + \sum_{i=1}^{k}\binom{n}{\frac{n+t-1}{2}+i} - \binom{n-t}{\frac{n-t-1}{2}+k}.\]
\end{conjecture}

Frankl~\cite{F} more modestly conjectures the following (Frankl's conjecture is formulated for $s$-union families rather than $t$-intersecting families, but our formulation is equivalent to Frankl's after taking complements). 

\begin{conjecture}[Frankl]\label{frankloddcase}
Let $g(n, t, k):= \max\{|\cG| - |\Delta_{\frac{n-t+1}{2} - k}(\cG)|: \hspace{2mm} \cG \subset \binom{[n]}{\frac{n-t+1}{2}}$ { is inter\-sec\-ting}$\}$. Then, if $n+t$ is odd and $\cF$ is a $t$-intersecting, $k$-Sperner family, then 
\[|\cF| \le g(n, t, k) + \sum_{i=1}^{k}\binom{n}{\frac{n+t-1}{2}+i}.\]
\end{conjecture}

Note that Conjecture~\ref{oddcase} can be interpreted as a strengthening of Conjecture~\ref{frankloddcase}, in that additionally there is a conjecture for the value of the function $g(n, t, k)$ for sufficiently large $n$. The connection may be made more apparent by noting that, after taking complements, we may equivalently define $g(n, t, k):= \max\{|\cG| - |\nabla_{\frac{n+t-1}{2}+k}(\cG)|: \hspace{2mm} \cG \subset \binom{[n]}{\frac{n+t-1}{2}} \text{ is $t$-intersecting}\}.$

\subsection{New Results} Let us mention that Frankl proved Conjecture \ref{evencase} when $t\ge n-O(\sqrt{n})$. We settle  Conjecture \ref{evencase} if $t$ is fixed and $n$ is sufficiently large. 

\begin{theorem}\label{main}
Let $t$ and $k$ be positive integers, and suppose that $n+t$ is even with $t\le n$, and $n$ is large enough. If $\cF\subseteq 2^{[n]}$ is a $t$-intersecting $k$-Sperner family, then 
$$|\cF|\le \binom{[n]}{\frac{n+t}{2}}+\ldots+\binom{[n]}{\frac{n+t}{2}+k-1}.$$
\end{theorem}

\section{Proof of Theorem~\ref{main}}

\subsection{Main ideas} The proof has two parts. In the first part we compress  $\cF$,  a $t$-intersecting, $k$-Sperner family, into the layers of the Boolean lattice containing the sets of sizes $\frac{n+t}{2}-k+1,\ldots,\frac{n+t}{2}+2k-2$. 
This part of the proof also works when $n+t$ is odd. In the second part, we use Katona's circle method, i.e., for every cyclic permutation $\sigma$ of $[n]$ we define $\cF_\sigma$ to be the collection of sets from $\cF$, whose elements are consecutive on $\sigma$, the so-called \textit{intervals}. For every $\sigma$, we show that for an appropriate weight function $w$, the total weight $w(\cF_\sigma)$ is maximized when $\cF_\sigma$ contains 
all intervals of size $r$ for every $\frac{n+t}{2}\le r\le \frac{n+t}{2}+k-1$. Then we deduce the general problem to this weighted version of the problem on the cycle.

\subsection{Compression Argument}

We recall the well-known Katona shadow $t$-intersection theorem~\cite{K}. 
\begin{theorem}[Katona shadow $t$-intersection theorem]\label{thm:kashadowtintersectionthm}
Let $\cF$ be an $r$-uniform, $t$-intersecting family. Then, for $r-t\le \ell\le r$, 
\[|\Delta_{\ell}(\cF)|\ge \frac{\binom{2r-t}{\ell}}{\binom{2r-t}{r}}|\cF|.\]
\end{theorem}

We prove a lemma about the $(i+1)$-shade of $\cF_i$ for $i\le \lfloor{\frac{n+t-1}{2}\rfloor}$.

\begin{lemma}[\cite{M}]\label{i+1shadelemgen} 
For $i\le \lfloor{\frac{n+t-1}{2}\rfloor}$, if $\cF \subseteq 2^{[n]}$ is $t$-intersecting, then we have 
\[|\nabla_{i+1}(\cF_i)| \ge |\cF_i|.\]
\end{lemma}

\begin{proof}[Proof]
Define the family of complements $\cF_i^C:= \{F^{C}: \hspace{1mm} F\in \cF_i\}$. Since $\cF_i$ is $t$-intersecting, $\cF_i^C$ is $(n+t-2i)$-intersecting. Since $i\le \lfloor{\frac{n+t-1}{2}\rfloor}$, we have $n+t-2i \ge 1$, so Theorem \ref{thm:kashadowtintersectionthm} can be applied to $\cF_i^C$ with $r:= n - i$, $t:= n+t-2i$, and $\ell:= n-i-1$, yielding 
\[
|\Delta_{n-i-1}(\cF_{i}^C)| \ge \frac{\binom{2(n-i) - (n+t-2i)}{n-i-1}}{\binom{2(n-i)-(n+t-2i)}{n-i}}|\cF_i^C|
					=\frac{\binom{n-t}{n-i-1}}{\binom{n-t}{n-i}}|\cF_i^C| = \frac{n-i}{i-t+1}|\cF_i|
					\ge |\cF_i|. 
\]
Since $|\nabla_{i+1}(\cF_i)| = |\Delta_{n-i-1}(\cF_{i}^C)|$, the desired result follows. 
\end{proof}

\begin{lemma}\label{tintk+1chfreven}
Let $\cF\subseteq 2^{[n]}$ be a $t$-intersecting and $k$-Sperner family, where   $n+t$ is even. Then there exists a $t$-intersecting $k$-Sperner family $\cG\subseteq 2^{[n]}$ with $|\cG| \ge |\cF|$ and $\min\{|G|: G\in \cG\} \ge \frac{n+t}{2} - (k-1)$. 
\end{lemma}

\begin{proof}
Recall that $f_i:=|\cF_i|$. 
Assume that there is $i < \frac{n+t}{2} - (k-1)$ such that $f_i > 0$ and $f_j = 0$ for every $j < i$. We show that there is a $t$-intersecting $k$-Sperner family $\cF'$ with $|\cF'| \ge |\cF|$ and $|F'| \ge i+1$ for every set $F'\in \cF'$. We show the existence of such an $\cF'$ by using a compression operation. 

We define a series of auxiliary families $\cH_j$ for $j\ge i$  as follows: $\cH_i := \cF_i$ and $\cH_j := \nabla_j(\cH_{j-1}) \cap \cF_{j}$ for $j > i$. The compression operation is as follows: we compress the sets in $\cH_i$ onto $\nabla_{i+1}(\cH_i)$. If $\cH_{i+1} = \emptyset$, we stop. Otherwise, we think of the sets of $\cH_{i+1}$ as appearing with multiplicity two in the newly constructed intermediate family. We compress one of the copies of each set in $\cH_{i+1}$ onto its $(i+2)$-shade $\nabla_{i+2}(\cH_{i+1})$, and leave the other copy on the $(i+1)$-layer. If $\cH_{i+2} \neq\emptyset$, then we repeat this compression process. We do the same for every $j \ge i$ as long as $\cH_{j} \neq \emptyset$. This compression process must terminate, since $\cH_{i+k} =\emptyset$, as otherwise there would be a $(k+1)$-chain in $\cF$. Call the family obtained after performing this series of compressions $\cF'$. 

In each step we added elements to the sets, hence  $\cF'$ will be $t$-intersecting. By Lemma~\ref{i+1shadelemgen}, we have $|\nabla_{j+1}(\cH_{j})| \ge |\cH_j|$ for $j \le \frac{n+t}{2} - 1$, so $|\cF'| \ge |\cF|$. It remains to be shown that $\cF'$ contains no $(k+1)$-chains. 

Let $A_0 \subset A_1 \subset \ldots \subset A_k$ be a $(k+1)$-chain in $\cF'$ with $|A_0| = j$, where $i+1\le j \le i+k$. Note that $|A_1| \ge j+1$, \ldots, $|A_{i+k-j}| \ge j+(i+k-j)=i+k$, so $A_{i+k-j+1}, \ldots, A_k \in \cF$. 

If all of the sets $A_0, A_1, \ldots, A_{i+k-j}$ were contained in $\cF$, then the $(k+1)$-chain $A_0 \subset \ldots \subset A_k$ would have already been in $\cF$. Otherwise, pick the largest $m$ such that $A_m \notin \cF$, and assume that $|A_m| = \ell$, so that $\ell \ge j+m$. By construction, there must be a chain $B_0 \subset B_1 \subset \ldots \subset B_{\ell-i-1} \subset A_m$ with $B_r \in \cF_{i+r}$ for $0\le r \le \ell-i-1$. Now the chain $B_0 \subset B_1 \subset \ldots \subset B_{\ell-i-1} \subset A_{m+1} \subset \ldots \subset A_k$ is contained in $\cF$, and it has size $\ell - i + k - m \ge j-i+k \ge k+1$, which is a contradiction.
\end{proof}

{\bf Remark.} The bottleneck of the proof is that we need to do the upshifting operation $k-1$ times, and we need the shade to be expanding, i.e., that is the reason that we require $i < \frac{n+t}{2} - (k-1)$. 
Observe that the parity of $n+t$ was not considered, so the same proof works when $n+t$ is odd. 

\begin{lemma}\label{easyshadow}
Let $\cF\subseteq 2^{[n]}$ be a Sperner family with $m:=\min\{|F|:F\in \cF\}>n/2$. Then for every $\lfloor n/2\rfloor \le j\le m$, we have $|\Delta_j(\cF)|\ge |\cF|$. \qed
\end{lemma}

Lemma~\ref{easyshadow} follows from a simple double-counting argument that was already used by Sperner~\cite{S} in his original proof.

\begin{lemma}\label{downshift}
If $\cF\subseteq 2^{[n]}$ is a $t$-intersecting $k$-Sperner family with $\min \{|F|:F\in \cF\}=\frac{n+t}{2}-c$, then there exists a $t$-intersecting $k$-Sperner family $\cF'\subseteq 2^{[n]}$ with $|\cF|\le |\cF'|$, and 
$$ \min\{|F|:F\in \cF\}=\min\{|F'|:F'\in \cF'\} \quad {\rm and }\quad \max \{|F'|:F'\in \cF'\}\le\frac{n+t}{2}+c+k-1.$$
\end{lemma}

\begin{proof}
We first partition $\cF$ into $\cF^1,\cF^2,\dots, \cF^k$ by letting $\cF^1$ consist of all minimal sets of $\cF$ and once $\cF^1, \ldots, \cF^j$ are defined, then let $\cF^{j+1}$ consist of all the minimal sets of $\cF\setminus \cup_{i=1}^j\cF^j$. Then for every $1\le j\le k$, we partition $\cF^j$ into $\cF^{j>}\cup \cF^{j\le}$ with 
$$\cF^{j>}=\{F\in \cF^j:|F|>\frac{n+t}{2}+c+j-1\}  \quad {\rm and }\quad \cF^{j\le}=\{F\in \cF^j:|F|\le\frac{n+t}{2}+c+j-1\}.$$ 
We define $\cF'^j:=\cF^{j\le}\cup \Delta_{\frac{n+t}{2}+c+j-1}(\cF^{j>})$.

Clearly, all the $\cF'^j$s are antichains. By Lemma \ref{easyshadow}, we have $|\cF^j|\le |\cF'^j|$ for all $1\le j\le k$ and thus for $\cF':=\cup_{j=1}^k\cF'^j$, we have $|\cF|\le |\cF'|$. 
Observe that  $\cF'$ is $k$-Sperner as it is the union of $k$ antichains. Additionally, $\cF'$ contains no set twice, since if $F$ were obtained after down-shifting of some $F\cup \{x\}$, then $F\in \cF$ would also have been down-shifted. Finally, $\cF'$ is $t$-intersecting as all sets in $\cF'\setminus\cF$ have size at least $\frac{n+t}{2}+c$ and all sets in $\cF\cap \cF'$ have size at least $\frac{n+t}{2}-c$.
\end{proof}

Observe that starting with an arbitrary $t$-intersecting $k$-Sperner family $\cF$, after applying Lemma \ref{tintk+1chfreven} we obtain another one $\cF'$ with $|\cF|\le |\cF'|$ and $\min \{|F|:F\in \cF\}\ge \frac{n+t}{2}-k+1$. Then applying Lemma \ref{downshift} with $c=\frac{n+t}{2}-\min\{|F|:F\in \cF'\}$, we obtain a $t$-intersecting $k$-Sperner family $\cF''$ with $|\cF|\le |\cF'|\le |\cF''|$ and $\min\{|F|: F\in \cF\}=\frac{n+t}{2}-m$ for some $0\le m\le k-1$ and $\max\{|F|:F\in \cF''\}\le \frac{n+t}{2}+k-1+m$. Therefore, in the next subsection, in the rest of the proof of Theorem \ref{main}, we will assume that $\cF$ has this property.

\subsection{Proof of Theorem \ref{main}}
Let $\sigma$ be a cyclic permutation of $[n]$ and $\cF_\sigma$ be the subfamily of those sets in $\cF$ that form an interval in $\sigma$. Note that there are $(n-1)!$ choices for $\sigma$. For a set $G$, let $w(G)=\binom{n}{|G|}$ and $w(\cG)=\sum_{G\in \cG}w(G)$.  We define $m$  as $m:=\frac{n+t}{2}- \min\{|F|: F\in \cF''\}$. By the above discussions, we have $0\le m\le k-1$. If $m=0$ then $\cF$ has the required structure, hence we assume $m>0$. The aim of this subsection is to prove the following lemma.

\begin{lemma}\label{cycle+}
Suppose $n+t$ is even with $t\le n$ and $n$ is large enough. For every cyclic permutation $\sigma$ and $t$-intersecting $k$-Sperner family $\cF\subseteq \bigcup_{i=\frac{n+t}{2}-m}^{\frac{n+t}{2}+k-1+m}\binom{[n]}{i}$, we have $w(\cF_\sigma)\le n\sum_{i=0}^{k-1}\binom{n}{\frac{n+t}{2}+i}$.
\end{lemma}

Before continuing, let us show how Lemma \ref{cycle+} implies Theorem \ref{main}.

\begin{proof}[Proof of Theorem \ref{main} using Lemma \ref{cycle+}]
As mentioned in the last paragraph of the previous subsection, by Theorem \ref{tintk+1chfreven} and Lemma \ref{downshift}, we can assume that $\cF\subseteq \bigcup_{i=\frac{n+t}{2}-m}^{\frac{n+t}{2}+k-1+m}\binom{[n]}{i}$ holds. Then using Lemma \ref{cycle+} we have:
\[
\sum_{\sigma}\sum_{F\in \cF_\sigma}w(F)\le (n-1)!\cdot n\sum_{i=0}^{k-1}\binom{n}{\frac{n+t}{2}+i}=n!\cdot\sum_{i=0}^{k-1}\binom{n}{\frac{n+t}{2}+i}.
\]
From the other side,
\[
\sum_{\sigma}\sum_{F\in \cF_\sigma}w(F)=\sum_{F\in \cF}|F|!(n-|F|)!\binom{n}{|F|}=n!|\cF|,
\]
which implies the required upper bound on $|\cF|$.
\end{proof}

In order to prove Lemma \ref{cycle+}, we need some preparation. Let us fix a cyclic permutation $\sigma$ of $[n]$. We partition all intervals, i.e., sets of consecutive elements of $[n]$ with respect to $\sigma$, into $n$ chains: the $h$-th chain $C_h$ consists of $\{\sigma(h)\},\{\sigma(h),\sigma(h+1)\},\dots,[n]\setminus\{\sigma(h-1)\}$ and we let $\cC_\sigma=\{C_h:h\in[n]\}$. A family $\cG$ of intervals is \textit{$\sigma$-$k$-Sperner $t$-intersecting} if it is $t$-intersecting and for every $C\in\cC_\sigma$ we have $|\cG\cap C|\le k$. Such a family is \textit{consecutive} if for every $C\in\cC_\sigma$ the chain $C\cap \cG$ consists of consecutive intervals, and \textit{full consecutive} if further $|C\cap \cG|=k$ holds for every $C\in\cC_\sigma$. Clearly, if $\cG$ is a $k$-Sperner $t$-intersecting family of intervals of $[n]$ with respect to $\sigma$, then $\cG$ is $\sigma$-$k$-Sperner $t$-intersecting, and if $\cF\subseteq 2^{[n]}$ is $t$-intersecting $k$-Sperner, then $\cF_\sigma$ is $\sigma$-$k$-Sperner $t$-intersecting for any $\sigma$. As the $t$-intersection property depends only on the smallest intervals $G_h\in C_h$, one can replace any $G\in \cG\cap C_h$ by any $G'\in C_h\setminus \cG$ with $|G'|>|G_h|$ to obtain another $\sigma$-$k$-Sperner $t$-intersecting family. So if $G'_h$ is the maximum interval of $C_h\cap \cG$ and $G$ is an interval from $C_h\setminus \cG$ with $|G_h|<|G|<|G'_h|$, then we can proceed as follows: if $|G|\ge n/2$, then we replace $G'_h$ by $G$, while if $|G|<n/2$, then we replace $G_h$ by $G$ to obtain a family $\cG'$. By our choice, we have $w(\cG)<w(\cG')$. As the difference of the sizes of maximum and minimum intervals of $C_h$ in the family strictly decreased, after a finite number of replacements, we obtain a consecutive $\sigma$-$k$-Sperner $t$-intersecting family. Note that during this process, we do not care if we created a $(k+1)$-chain which is not in one of the $C_h$. 

Finally, we can extend any consecutive  $\sigma$-$k$-Sperner $t$-intersecting family to a full one. 
More generally, the following holds.

\begin{observation}
If $\cG$ is a $k$-Sperner $t$-intersecting family of intervals with respect to $\sigma$ such that every interval has size between $\frac{n+t}{2}-m$ and $\frac{n+t}{2}+k-1+m$ for some $0\le m\le k-1$, then there exists a full consecutive $\sigma$-$k$-Sperner $t$-intersecting family $\cG'$ with $w(\cG)\le w(\cG')$.
\end{observation}

\begin{proof}
By the argument above, we can obtain a consecutive $\sigma$-$k$-Sperner $t$-intersecting family $\cG^*$. If for a chain $C_h$, we have $|\cG^*\cap C_h|<k$, then we add the interval of $C_h$ to $\cG^*$ that is one larger than the maximum interval in $\cG^*\cap C_h$. We could only get into trouble if for some $h$ the smallest interval $G_h$ of $C_h\cap \cG^*$ is strictly larger than $\frac{n+t}{2}+m$, but then we can add $G\in C_h$ with $|G|=\frac{n+t}{2}+m$ to $\cG^*$ without violating the $t$-intersecting property as such $G$ $t$-intersects all other intervals in $\cG^*$ because of the size restrictions. 
\end{proof}

To prove Lemma \ref{cycle+}, it is sufficient to show the following statement.

\begin{lemma}\label{generalcycle}
Suppose $n+t$ is even and $n$ is large enough. Let $\cG$ be a full consecutive  $\sigma$-$k$-Sperner $t$-intersecting family of intervals on a cycle of length $n$ such that $\min\{|G|:G\in \cG\}=\frac{n+t}{2}-m$ and  $\max\{|G|:G\in \cG\}\le \frac{n+t}{2}+k-1+m$ for some $0\le m\le k-1$. Then $w(\cG)\le n\sum_{i=0}^{k-1}\binom{n}{\frac{n+t}{2}+i}$ holds.
\end{lemma}

The following Fact is easy to see, and was known previously, we omit its proof.

\begin{fact}\label{shad}
Let $\cG$ be a family of intervals on a cycle of length $n$. If $\cG$ consists of $i$-intervals for some $2\le i\le n-1$, then $|\Delta(\cG)|\ge |\cG|$. 

\end{fact}

The next Fact is standard, we include its proof, as it is a bit technical to see it instantly.

\begin{fact}\label{binom}
For every fixed $a<b$ with $0<b$ there exists $n_0=n_0(a,b)$ such that if $n\ge n_0$, then we have $\binom{n}{\lfloor n/2\rfloor +a}+\binom{n}{\lfloor n/2\rfloor +b}\le \binom{n}{\lfloor n/2\rfloor +a+1}+\binom{n}{\lfloor n/2\rfloor +b-1}$.
\end{fact}

\begin{proof}
If $a$ is negative, then $\binom{n}{\lfloor\frac{n}{2}\rfloor +a}<\binom{n}{\lfloor\frac{n}{2}\rfloor +a+1}$ and $\binom{n}{\lfloor\frac{n}{2}\rfloor +b}\ge \binom{n}{\lfloor\frac{n}{2}\rfloor +b-1}$ so the statement of the Fact holds for all values of $n$.

Hence, we can assume $0\le a<b$. If $b=a+1$, then clearly equality holds.
If $b>a+1$, then dividing by $n!$ and multiplying by $(\lfloor n/2\rfloor+a+1)!\cdot (\lceil n/2\rceil-a)!\cdot (\lfloor n/2\rfloor+b)!\cdot (\lceil n/2\rceil-b+1)!$, the desired inequality is equivalent to
\[
(\lfloor n/2\rfloor+a+1)\cdot (\lfloor n/2\rfloor+b)!\cdot (\lceil n/2\rceil-b+1)!+
(\lfloor n/2\rfloor+a+1)!\cdot (\lceil n/2\rceil-a)!\cdot (\lceil n/2\rceil-b+1)\le 
\]
\[
 (\lceil n/2\rceil-a)\cdot (\lfloor n/2\rfloor+b)!\cdot (\lceil n/2\rceil-b+1)!+ (\lfloor n/2\rfloor+a+1)!\cdot (\lceil n/2\rceil-a)!\cdot (\lfloor n/2\rfloor+b).
\]
Rearranging gives
\[
(2a+1)\cdot (\lfloor n/2\rfloor+b)!\cdot (\lceil n/2\rceil-b+1)!\le (2b-1)\cdot (\lfloor n/2\rfloor+a+1)!\cdot (\lceil n/2\rceil-a)!,
\]
which is equivalent to
\[
\frac{2a+1}{2b-1}\le \frac{(\lceil n/2\rceil-a)\cdot \ldots \cdot (\lceil n/2\rceil-b+2)}{(\lfloor n/2\rfloor +b)\cdot \ldots \cdot (\lfloor n/2\rfloor +a+2)}.
\]
The left hand side is a fixed rational number smaller than 1, while the right hand side tends to one as $n$ tends to infinity.
\end{proof}

The next simple observation is going to be the core of our argument. For a cyclic permutation $\sigma$ and an interval $G$ define $\overline{G}^t$ as the complement of $G$ together with the (counterclockwise) leftmost $\lfloor \frac{t}{2}\rfloor$ and rightmost $\lceil 
\frac{t}{2}\rceil$ elements of $G$ with respect to $\sigma$. For a family $\cG$ of intervals, let $\overline{\cG}^t=\{\overline{G}^t :G\in\cG\}$.

\begin{lemma}\label{complement+}
Suppose $n+t$ is even, $\sigma$ is a cyclic permutation of $[n]$. If $\cG$ is a full consecutive  $\sigma$-$k$-Sperner $t$-intersecting family with interval sizes between $\frac{n+t}{2}-m$ and $\frac{n+t}{2}+k-1+m$ for some $0\le m\le k-1$, then for any $G \in \cG$ no proper subinterval $H$ of $\overline{G}^t$ belongs to $\cG$. Moreover, if $|G|=\frac{n+t}{2}-m$, then $\overline{G}^t\in \cG$.
\end{lemma}

\begin{proof}
 Any proper subinterval of $\overline{G}^t$ intersects $G$ in less than $t$ elements, thus $\overline{G}^t$ cannot contain any interval from $\cG$.

As $|G|+|\overline{G}^t|=n+t$, if $|G|=\frac{n+t}{2}-m$, then $\overline{G}^t$ $t$-intersects every element of $\cG$.
 Also, if $\overline{G}^t\notin \cG$, then for the chain $C_h$ containing $\overline{G}^t$, we have $|\cG\cap C_h|<k$ as there are $k-1$ intervals in $C_h$ that are larger than $\overline{G}^t$. This contradicts the full consecutive $\sigma$-$k$-Sperner property.\end{proof}

The next lemma establishes some inequalities on the number of intervals that a full consecutive $\sigma$-k-Sperner $t$-intersecting family $\cG$ satisfying the assumptions of Lemma \ref{complement+} may contain. For $i=-m,-m+1,\dots,k+m-1$, let $\cG_i$ be the family of intervals of length $\frac{n+t}{2}+i$ in $\cG$ and let $ g_i$ denote the size of $\cG_i$.

\begin{lemma}\label{inequalities}
Suppose $n+t$ is even, $m<k$ and $n$ is large enough. Let $\cG$ be a full consecutive  $\sigma$-$k$-Sperner $t$-intersecting family of intervals on a cycle of length $n$ such that $\min\{|G|:G\in \cG\}=\frac{n+t}{2}-m$ and  $\max\{|G|:G\in \cG\}\le \frac{n+t}{2}+k-1+m$. Then we have the following inequalities:
\begin{enumerate}
    \item
    $g_{-j-1}+g_j\le n$ for all $0\le j\le m-1$ satisfying  $j<k-m$.
    \item
    $\sum_{i=-(j+1)}^jg_i\le (j+1)n-\sum_{i=k-j}^{m}g_{-i}$ for all $0\le j\le m-1$ such that $j\ge k-m$.
    \item 
    $\sum_{i=-m}^{k-j}g_i\le (k-j+1)n-\sum_{i=j}^mg_{-i}  $ for all $1\le j\le m$ such that $j<k-m$. 
    \item
    $\sum_{i=-m}^{k-1+j-1}g_i\le kn-\sum_{i=-m}^{-j}g_i$ for all $1\le j\le m$.
\end{enumerate}
\end{lemma}

\begin{proof}
To prove (1), note that Lemma \ref{complement+} applied to $\cG_{-(j+1)}$ implies that $\Delta(\overline{\cG_{-(j+1)}}^t)$ is disjoint from $\cG$. Fact \ref{shad} implies that $|\Delta(\overline{\cG_{-(j+1)}}^t)|\ge |\overline{\cG_{-(j+1)}}^t| = g_{-j-1}$, and since $\Delta(\overline{\cG_{-(j+1)}}^t)$ is a family of $(\frac{n+t}{2}+j)$-intervals, it follows that $g_j + g_{-j-1} \le g_j + |\Delta(\overline{\cG_{-(j+1)}}^t)| \le n$. 

The proofs of (2) and (3) are similar. We define families $\cH_1,\cH_2$ of missing intervals, i.e. that are not members of $\cG$, as follows: 
\[
\cH_1^j=\{H\notin \cG: \frac{n+t}{2}\le |H| \le \frac{n+t}{2}+j, ~\not\exists G\in \cG \hskip 0.15truecm H \supset G\}, 
\]
\[\cH_2^j=\{H\notin \cG: \frac{n+t}{2}\le |H| \le \frac{n+t}{2}+j, ~\not\exists G\in \cG \hskip 0.15truecm G \supset H\}.
\]
Observe that by definition and by the full consecutive property, we have $\cH_1\cap \cH_2=\emptyset$.

To prove (2), we consider $\cH_1^j$ and $\cH_2^j$. First, as in the proof of (1), for any $1\le i\le j+1$, $\Delta(\overline{\cG_{-i}}^t)$ is disjoint from $\cG$ by Lemma \ref{complement+} and, by Fact \ref{shad}, has size at least $g_{-i}$. Note that all these missing intervals (missing from $\cG$) are \textit{only below} intervals of $\cG$, so $\cH_1^j\supseteq \bigcup_{i=1}^{j+1}\Delta(\overline{G_{-i}}^t)$.
On the other hand, if $G\in \cG_{-i}$ with $k-j\le i\le m$, then, as $\cG$ is consecutive, the chain $C_h\in \cC_\sigma$ that contains $G$ misses all intervals that are exactly $k$ larger than $|G|$, i.e. of size  $\frac{n+t}{2}+k-i$. 
 Thus we obtain that $\cH_2^j$ contains at least $\sum_{i=k-j}^mg_{-i}$ missing intervals each of which are \textit{only above} some intervals of $\cG$. This means that $\cup_{i=0}^j\cG_i$, $\cH_1^j$, $\cH_2^j$ are pairwise disjoint, have sizes $\sum_{i=0}^jg_i$, $\sum_{-(j+1)}^{-1}g_i$, and $\sum_{i=k-j}^mg_{-i}$, and contain intervals of sizes between $\frac{n+t}{2}$ and $\frac{n+t}{2}+j$. There are $(j+1)n$ such intervals, therefore $\sum_{i=0}^jg_i+\sum_{-(j+1)}^{-1}g_i+\sum_{i=k-j}^mg_{-i}\le (j+1)n$ holds. Merging the first two terms and rearranging yields (2).

To prove (3), we consider $\cH_1^{k-j},\cH_2^{k-j}$. As $j<k-m$, this time $\cup_{i=1}^m\Delta(\overline{\cG_{-i}}^t)$ belongs to $\cH_1^{k-j}$, and by Lemma \ref{complement+} and Fact~\ref{shad}, $\cH_1^{k-j}$ has size at least $\sum_{i=1}^mg_{-i}$. For any $G\in \cG_{-i}$ with $i\ge j$, the intervals $G'$ of $C_h$ with $G\in C_h$ and $|G'|-|G|\ge k$ are missing by the consecutive property of $\cG$. There are $i-j+1$ of such missing intervals. We obtain that $\cH_2$ contains at least $\sum_{i=j}^mg_{-i}$ missing intervals. Again, $\cup_{i=0}^{k-j}\cG_i$, $\cH_1^{k-j}$, and $\cH_2^{k-j}$ are pairwise disjoint, so the sum of their sizes is at most $(k-j+1)n$. After rearrangement, this yields (3). 

Finally, to see (4) observe first that as for a full consecutive $\sigma$-$k$-Sperner $t$-intersecting family, we have $|C_h\cap \cG|=k$ for all $h$, we have $\sum_{i=-m}^{k-1+m}g_i=|\cG|=kn$. So the statement of (4) is equivalent to the statement that the number of intervals in $\cG$ of size at least $\frac{n+t}{2}+k-1+j$ is at least $\sum_{i=-m}^{-j}g_i$. Again, we apply Lemma \ref{complement+} and observe that intervals $G$ of $\overline{\cG_{-i}}^t$ do not strictly contain any interval of $\cG$. Therefore, if $i\ge j$, then the chain $C_h$ containing $G$ has at least $i-j+1$ intervals from $\cG$ of size at least $\frac{n+t}{2}+k+j-1$. Counting all these, $\cG$ contains at least $\sum_{i=-m}^{-j}g_i$ intervals of size at least $\frac{n+t}{2}+k+j-1$ as desired.
\end{proof}

We are now ready to prove Lemma \ref{generalcycle}. As the proof involves lots of formulas, we sketch the main idea. As mentioned in the last paragraph of the proof of Lemma \ref{inequalities}, the size of a full consecutive $\sigma$-$k$-Sperner $t$-intersecting family $\cG$ is $kn$, so its weight $w(\cG)=\sum_{G\in \cG}w(G)=\sum_{G\in\cG}\binom{n}{|G|}=\sum_{i=-m}^{k+m-1}g_i\binom{n}{\frac{n+t}{2}+i}$ is a sum of $kn$ binomial coefficients. If all $g_i$s are 0 whenever $i$ is negative, then we are done as $w(|G|)$ is monotone decreasing in $|G|$ if $|G|>\frac{n}{2}$, and $g_i\le n$ for all $i$. If there exists $i<0$ with $g_i>0$, then we plan to apply Fact \ref{binom} to obtain another set of coefficients $g'_i$ such that $\sum_i g'_i=kn$, $\sum_ig_i\binom{n}{\frac{n+t}{2}+i}\le \sum_{i}g'_i\binom{n}{\frac{n+t}{2}+i}$, and $\sum_{i=0}^jg'_i\le (j+1)n$ hold for all $j=0,1,\dots, k-1$. When applying Fact \ref{binom}, we will match $g_{-i}$ with $g_{k-1+i}$ for all $i=1,2,\dots,m$, therefore we will first have to make sure that $g_{-i}\le g_{k-1+i}$ and then we can apply Fact~\ref{binom}.

\begin{proof}[Proof of Lemma \ref{generalcycle}]
Note that if $m=0$ then we are done, hence we assume $m\ge 1$. 
Let us introduce the coefficients $g'_i$:
\begin{itemize}
    \item 
    for $-m\le i\le k-1$, let $g'_i=g_i$,
    \item
    for $i=2,3,\dots,m$ let $g'_{k-1+i}=g_{-i}$,
    \item
    let $g'_k=\sum_{i=k}^{k+m-1}g_i-\sum_{i=2}^mg_{-i}$.
\end{itemize}
By the definition of $g'_k$ and the fact that $|\cG|=\sum_{i=-m}^{k+m-1}g_i= kn$, we have $\sum_{i=-m}^{k+m-1}g'_i=kn$ and $\sum_{i=k}^{k+m-1}g'_i =  \sum_{i=k}^{k+m-1}g_i$.

Lemma \ref{inequalities} (4) states $\sum_{i=-m}^{k-1+j-1}g_i\le kn-\sum_{i=-m}^{-j}g_i$. Plugging in $kn=\sum_{i=-m}^{k-1+m}g_i$ and rearranging yields
\begin{equation}\label{4equiv}
  \sum_{i=-m}^{-j}g_i\le \sum_{i=k-1+j}^{k-1+m}g_i \tag{*}
\end{equation}

Applying (\ref{4equiv}) with $j=1$, we obtain 
\begin{equation}\label{g_minus}g'_k=\sum_{i=k}^{k+m-1}g'_i -\sum_{i=k+1}^{k-1+m}g'_i
= \sum_{i=k}^{k+m-1}g_i   -\sum_{i=-m}^{-2}g'_i \ge \sum_{i=-m}^{-1}g_i-\sum_{i=-m}^{-2}g_i=g_{-1}.\tag{**} 
\end{equation}
We would like to compare $\sum_{i=k}^{k+m-1}g'_i\binom{n}{\frac{n+t}{2}+i}$ to $\sum_{i=k}^{k+m-1}g_i\binom{n}{\frac{n+t}{2}+i}$. As mentioned above, $A:=\sum_{i=k}^{k+m-1}g'_i=\sum_{i=k}^{k+m-1}g_i$. Also, (\ref{4equiv}) and $g'_{k-1+i}=g_{-i}$ for all $i=2,3,\dots,m$ imply 
\begin{equation}\label{dia}
\sum_{i=k-1+j}^{k-1+m}g'_i= \sum_{i=-m}^{-j}g_i\le \sum_{i=k-1+j}^{k-1+m}g_i\tag{$\diamond$}
\end{equation}
for all $j=2,3,\dots,m$. Therefore, we can apply the following general statement that can be easily seen by induction: $a_1,a_2,\dots,a_n$, $b_1,b_2,\dots,b_n$ and $d_1\ge d_2\ge \dots \ge d_n$ are all non-negative integers with $\sum_{i=j}^na_i\le \sum_{i=j}^nb_i$ for all $j=2,2,\dots,n$ and $\sum_{i=1}^na_i= \sum_{i=1}^nb_i$. Then $\sum_{i=1}^na_id_i\ge \sum_{i=1}^nb_id_i$.
Plugging in $n:=m$, $a_i:=g'_i$, $b_i:=g_i$, and $d_i:=\binom{n}{\frac{n+t}{2}+i}$, we obtain

$$\sum_{i=k}^{k+m-1}g'_i\binom{n}{\frac{n+t}{2}+i}\ge \sum_{i=k}^{k+m-1}g_i\binom{n}{\frac{n+t}{2}+i}.$$ Thus \begin{equation}\label{begin}
    w(\cG)=\sum_{i=-m}^{k+m-1}g_i\binom{n}{\frac{n+t}{2}+i}\le\sum_{i=-m}^{k+m-1}g'_i\binom{n}{\frac{n+t}{2}+i}.\tag{***}
\end{equation}
Now for every $1\le j\le m$, we apply~Fact \ref{binom} either $2j-1$ times if $\frac{n+t}{2}+j-1\le \frac{n+t}{2}+k-j$ or $k$ times if $\frac{n+t}{2}+j-1> \frac{n+t}{2}+k-j$  to obtain
\begin{equation}\label{push}
    g'_j\left(\binom{n}{\frac{n+t}{2}-j}+\binom{n}{\frac{n+t}{2}+k+j-1}\right)\le g'_j\left(\binom{n}{\frac{n+t}{2}+j-1}+\binom{n}{\frac{n+t}{2}+k-j}\right).\tag{****}
\end{equation}
Based on (\ref{push}), we want to give the weights of intervals of length $\frac{n+t}{2}-j$ to ``imaginary" intervals of length $\frac{n+t}{2}+j-1$ and those of length $\frac{n+t}{2}+k+j-1$ to those of length $\frac{n+t}{2}+k-j$. As $m\le k-1$, all imaginary intervals will have length between $\frac{n+t}{2}$ and $\frac{n+t}{2}+k-1$ (actually, $m\le k$ would suffice). Recall that $g'_j$ is defined for $-m\le j\le k+m-1$. Therefore, we introduce
\begin{eqnarray*}
g''_j=\left\{
\begin{array}{cc} 
g'_j+g'_{-(j+1)} & \textnormal{if}\ 0\le j\le m-1 ~\textnormal{and}\ \frac{n+t}{2}+j<\frac{n+t}{2}+k-m, \\
g'_j+g'_{-(j+1)}+g'_{-(k-j)} & \textnormal{if}\ 0\le j\le m-1 ~\textnormal{and}\ \frac{n+t}{2}+j\ge\frac{n+t}{2}+k-m,
\\
g'_j & \textnormal{if}\ j\ge m ~\text{and}\ \frac{n+t}{2}+j<\frac{n+t}{2}+k-m,
\\
g'_j+g'_{-(k-j)} & \textnormal{if}\ j\ge m ~\text{and}\ \frac{n+t}{2}+j\ge\frac{n+t}{2}+k-m,
\\
0 & \textnormal{if}\ j<0 ~\text{or}\ j>k\\
g'_k-g_{-1} & \textnormal{if}\ j=k.
\end{array}
\right.
\end{eqnarray*}
Observe that $g''_k=g'_k-g_{-1}\ge 0$, see~\eqref{g_minus}. Also, the values $g'_{-j}$ decreased to $0=g''_{-j}$ for $j=1,2,\dots,m$, but according to the first two cases of the definition of $g''$, the value of $g'_{j-1}$ was increased by $g_{-j}$ (and possibly something else). Also, the values of $g'_{k-1+j}=g'_{-j}$ were erased for $j=2,3,\dots,m$ and were given to $g'_{k-j}$ according to the second and fourth cases of the definition of $g''$. Finally, $g_{-1}=g'_{-1}$ from $g'_k$ was given to $g'_{k-1}$ according to the fourth case of the definition of $g''$. 
So $\sum_{i=0}^kg''_i=\sum_{i=-m}^{k-1+m}g'_i=kn$. 

Now, (\ref{push}) implies that (\ref{begin}) continues as 
\begin{equation}\label{continue}
    w(\cG)=\sum_{i=-m}^{k+m-1}g_i\binom{n}{\frac{n+t}{2}+i}\le\sum_{i=-m}^{k+m-1}g'_i\binom{n}{\frac{n+t}{2}+i}\le \sum_{i=-m}^{k+m-1}g''_i\binom{n}{\frac{n+t}{2}+i}.\tag{*****}
\end{equation}
We claim that for any $j=0,1,\dots,k-1$, we have $\sum_{i=0}^jg''_i\le (j+1)n$.

If $0\le j\le m-1$ such that $\frac{n+t}{2}+j<\frac{n+t}{2}+k-m$, then more is true: $$g''_j=g'_j+g'_{-(j+1)}=g_j+g_{-(j+1)}\le n$$ by Lemma \ref{inequalities} (1).
    
If $0\le j\le m-1$ such that $\frac{n+t}{2}+j\ge \frac{n+t}{2}+k-m$, then \begin{align*}
    \sum_{i=0}^jg''_i = & \sum_{i=0}^{k-1-m}g''_i+\sum_{i=k-m}^jg''_i \\
    = & \sum_{i=0}^{k-1-m}g_i+g_{-(i+1)}+\sum_{i=k-m}^jg_i+g_{-(i+1)}+g_{-(k-i)}=\sum_{-(j+1)}^jg_i+\sum_{i=k-j}^mg_{-i}\le (j+1)n
\end{align*} ensured by Lemma \ref{inequalities} (2).
    
If $m\le j<k-m$, then $g''_j=g'_j=g_j\le n$, so the inequality $\sum_{i=0}^jg''_i\le (j+1)n$ holds, as it holds in the previous two cases.
    
Finally, let us consider the case $j\ge \max\{m,\ge k-m\}$. If $j=k-1$, then there is nothing to prove as $\sum_{i=0}^{k-1}g''_i\le\sum_{i=0}^{k}g''_i=kn$. If $j\le k-2$, then we write $j=k-j^*$ with $2\le j^*\le m$ and obtain
\[
\sum_{i=0}^{k-j^*}g''_i=\sum_{i=-m}^{k-j^*}g_i+\sum_{i=j^*}^mg_{-i}\le (k-j^*+1)n
\]
by Lemma \ref{inequalities} (3).

To finish the proof of the lemma, observe that $g_i=0$ for $i<0$, $\sum_{i=0}^kg''_i= kn$, and $\sum_{i=0}^jg''_i\le (j+1)n$ imply that (\ref{continue}) can be continued as
\begin{equation*}
    w(\cG)=\sum_{i=-m}^{k+m-1}g_i\binom{n}{\frac{n+t}{2}+i}\le\sum_{i=-m}^{k+m-1}g'_i\binom{n}{\frac{n+t}{2}+i}\le \sum_{i=-m}^{k+m-1}g''_i\binom{n}{\frac{n+t}{2}+i}\le n\sum_{i=0}^{k-1}\binom{n}{\frac{n+t}{2}+i},
\end{equation*}
as claimed.
\end{proof}


{\bf Acknowledgements.} The authors are thankful to Adam Zs.~Wagner for fruitful discussions on the project.

\end{document}